\documentclass[leqno,12pt]{amsart}
\usepackage{amssymb} 
\theoremstyle{plain} 
\newtheorem{theorem}{\indent\sc Theorem}[section] 

\newtheorem{proposition}[theorem]{\indent\sc Proposition}

\theoremstyle{definition} 

\begin{document}

\title[Coupled Painlev\'e systems]{Coupled Painlev\'e systems with affine Weyl group symmetry of types $D_3^{(2)}$ and $D_5^{(2)}$ \\}
\author{Yusuke Sasano }

\renewcommand{\thefootnote}{\fnsymbol{footnote}}
\footnote[0]{2000\textit{ Mathematics Subjet Classification}.
34M55; 34M45; 58F05; 32S65.}

\keywords{ 
Affine Weyl group, birational symmetry, coupled Painlev\'e system.}
\maketitle

\begin{abstract}
In this paper, we find a two-parameter family of coupled Painlev\'e systems in dimension four with affine Weyl group symmetry of type $D_3^{(2)}$. We also find a four-parameter family of 2-coupled $D_3^{(2)}$-systems in dimension eight with affine Weyl group symmetry of type $D_5^{(2)}$.  We show that for each system, we give its symmetry and holomorphy conditions, respectively. These symmetries, holomorphy conditions and invariant divisors are new.
\end{abstract}

\section{Introduction}
In \cite{Sasa4,Sasa7,Sasa9}, we presented some types of coupled Painlev\'e systems with various affine Weyl group symmetries. In this paper, we find a 2-parameter family of coupled Painlev\'e systems in dimension four with affine Weyl group symmetry of type $D_3^{(2)}$ explicitly given by
\begin{align}
\begin{split}
\frac{dq_1}{dt}&=\frac{\partial H_{D_3^{(2)}}}{\partial p_1}, \quad \frac{dp_1}{dt}=-\frac{\partial H_{D_3^{(2)}}}{\partial q_1}, \quad \frac{dq_2}{dt}=\frac{\partial H_{D_3^{(2)}}}{\partial p_2}, \quad \frac{dp_2}{dt}=-\frac{\partial H_{D_3^{(2)}}}{\partial q_2}
\end{split}
\end{align}
with the polynomial Hamiltonian
\begin{align}\label{12}
\begin{split}
&H_{D_3^{(2)}}(q_1,p_1,q_2,p_2,t;\alpha_0,\alpha_1)\\
&=2H_{II}(q_1,p_1,t;\alpha_0)+H_{II}^{auto}(q_2,p_2;\alpha_1)+4p_1p_2-2q_1q_2p_2\\
&=2(q_1^2p_1+p_1^2+tp_1+\alpha_0 q_1)+q_2^2p_2-2p_2^2+\alpha_1 q_2+4p_1p_2-2q_1q_2p_2.
\end{split}
\end{align}
Here $q_1,p_1,q_2$ and $p_2$ denote unknown complex variables, and $\alpha_0,\alpha_1,\alpha_2$ are complex parameters satisfying the relation:
\begin{equation}\label{13}
\alpha_0+\alpha_1+\alpha_2=\frac{1}{2}.
\end{equation}
The symbols $H_{II}$ and $H_{II}^{auto}$ denote
\begin{align}
\begin{split}
H_{II}(x,y,t;\alpha_0)=&x^2y+y^2+ty+\alpha_0x\\
H_{II}^{auto}(z,w;\alpha_1)=&z^2w-2w^2+\alpha_1 z,
\end{split}
\end{align}
where $H_{II}$ denotes the second Painlev\'e Hamiltonian, and $H_{II}^{auto}$ denotes an autonomous version of the second Painlev\'e Hamiltonian. Of course, the system with the Hamiltonian $H_{II}^{auto}$ has itself as its first integral.

This is the first example which gave higher order Painlev\'e type systems of type $D_{3}^{(2)}$.

We remark that for this system we tried to seek its first integrals of polynomial type with respect to $q_1,p_1,q_2,p_2$. However, we can not find. Of course, the Hamiltonian $H_{D_3^{(2)}}$ is not the first integral.

We also find a 4-parameter family of 2-coupled $D_3^{(2)}$-systems in dimension eight with affine Weyl group symmetry of type $D_5^{(2)}$ explicitly given by
\begin{align}
\begin{split}
\frac{dq_1}{dt}&=\frac{\partial H_{D_5^{(2)}}}{\partial p_1}, \quad \frac{dp_1}{dt}=-\frac{\partial H_{D_5^{(2)}}}{\partial q_1}, \ldots ,\frac{dq_4}{dt}=\frac{\partial H_{D_5^{(2)}}}{\partial p_4}, \quad \frac{dp_4}{dt}=-\frac{\partial H_{D_5^{(2)}}}{\partial q_4}
\end{split}
\end{align}
with the polynomial Hamiltonian
\begin{align}\label{22}
\begin{split}
H_{D_5^{(2)}}=&H_{D_3^{(2)}}(q_1,p_1,q_2,p_2,t;\alpha_0,\alpha_1)+H_{D_3^{(2)}}(q_4,p_4,q_3,p_3,t;\alpha_4,\alpha_3)\\
&-\frac{3}{2}p_1^2-\frac{3}{2}p_4^2+3p_1p_4.
\end{split}
\end{align}
Here $q_1,p_1,q_2,p_2,q_3,p_3,q_4$ and $p_4$ denote unknown complex variables, and $\alpha_0,\alpha_1,\ldots,\alpha_4$ are complex parameters satisfying the relation:
\begin{equation}\label{23}
\alpha_0+\alpha_1+\alpha_2+\alpha_3+\alpha_4=1.
\end{equation}
We remark that for this system we tried to seek its first integrals of polynomial type with respect to $q_1,p_1,\ldots,q_4,p_4$. However, we can not find. Of course, the Hamiltonian $H_{D_5^{(2)}}$ is not the first integral.

This is the second example which gave higher order Painlev\'e type systems of type $D_{5}^{(2)}$.

We also remark that 2-coupled Painlev\'e III system in dimension four given in the paper \cite{Sasa8} admits the affine Weyl group symmetry of type $D_5^{(2)}$ as the group of its B{\"a}cklund transformations, whose generators $s_1,s_2,s_3$ are determined by the invariant divisors. However, the transformations $s_0,s_4$ do not satisfy so (see Theorem 4.1 in \cite{Sasa8}).

On the other hand, the system \eqref{21} admits the affine Weyl group symmetry of type $D_5^{(2)}$ as the group of its B{\"a}cklund transformations, whose generators $s_0,\ldots,s_4$ are determined by the invariant divisors \eqref{invariant} (see Section 3).

We show that for each system, we give its symmetry and holomorphy conditions, respectively. These B{\"a}cklund transformations of each system satisfy
\begin{equation}\label{universal}
s_i(g)=g+\frac{\alpha_i}{f_i}\{f_i,g\}+\frac{1}{2!} \left(\frac{\alpha_i}{f_i} \right)^2 \{f_i,\{f_i,g\} \}+\cdots \quad (g \in {\Bbb C}(t)[q_1,p_1,\ldots,q_4,p_4]),
\end{equation}
where poisson bracket $\{,\}$ satisfies the relations:
$$
\{p_1,q_1\}=\{p_2,q_2\}=\{p_3,q_3\}=\{p_4,q_4\}=1, \quad the \ others \ are \ 0.
$$
Since these B{\"a}cklund transformations have Lie theoretic origin, similarity reduction of a Drinfeld-Sokolov hierarchy admits such a B{\"a}cklund symmetry.

These symmetries, holomorphy conditions and invariant divisors are new.

\section{$D_3^{(2)}$ system} 
In this paper, we study a 2-parameter family of coupled Painlev\'e systems in dimension four with affine Weyl group symmetry of type $D_3^{(2)}$ explicitly given by
\begin{equation}\label{11}
  \left\{
  \begin{aligned}
   \frac{dq_1}{dt} &=\frac{\partial H_{D_3^{(2)}}}{\partial p_1}=2q_1^2+4p_1+2t+4p_2,\\
   \frac{dp_1}{dt} &=-\frac{\partial H_{D_3^{(2)}}}{\partial q_1}=-4q_1p_1+2q_2p_2-2\alpha_0,\\
   \frac{dq_2}{dt} &=\frac{\partial H_{D_3^{(2)}}}{\partial p_2}=q_2^2-4p_2+4p_1-2q_1q_2,\\
   \frac{dp_2}{dt} &=-\frac{\partial H_{D_3^{(2)}}}{\partial q_2}=-2q_2p_2+2q_1p_2-\alpha_1\\
   \end{aligned}
  \right. 
\end{equation}
with the polynomial Hamiltonian \eqref{12}.

\begin{theorem}
This system \eqref{11} admits the affine Weyl group symmetry of type $D_3^{(2)}$ as the group of its B{\"a}cklund transformations, whose generators are explicitly given as follows{\rm : \rm} with {\it the notation} $(*):=(q_1,p_1,q_2,p_2,t;\alpha_0,\alpha_1,\alpha_2),$
\begin{align}
\begin{split}
s_0:(*) \rightarrow &\left(q_1+\frac{4\alpha_0}{4p_1+q_2^2},p_1,q_2,p_2-\frac{2\alpha_0 q_2}{4p_1+q_2^2},t;-\alpha_0,\alpha_1+2\alpha_0,\alpha_2 \right),\\
s_1:(*) \rightarrow &\left(q_1,p_1,q_2+\frac{\alpha_1}{p_2},p_2,t;\alpha_0+\alpha_1,-\alpha_1,\alpha_2+\alpha_1 \right),\\
s_2:(*) \rightarrow &(q_1+\frac{4\alpha_2}{4p_1+8p_2+4q_1q_2-q_2^2+4t},p_1-\frac{4\alpha_2 q_2}{4p_1+8p_2+4q_1q_2-q_2^2+4t}\\
&-\frac{16{\alpha_2}^2}{(4p_1+8p_2+4q_1q_2-q_2^2+4t)^2},q_2+\frac{8\alpha_2}{4p_1+8p_2+4q_1q_2-q_2^2+4t},\\
&p_2-\frac{2\alpha_2(2q_1-q_2)}{4p_1+8p_2+4q_1q_2-q_2^2+4t},t;\alpha_0,\alpha_1+2\alpha_2,-\alpha_2).
\end{split}
\end{align}
\end{theorem}
Since these B{\"a}cklund transformations have Lie theoretic origin, similarity reduction of a Drinfeld-Sokolov hierarchy admits such a B{\"a}cklund symmetry.

\begin{proposition}
This system has the following invariant divisors\rm{:\rm}
\begin{center}\label{19}
\begin{tabular}{|c|c|c|} \hline
parameter's relation & $f_i$ \\ \hline
$\alpha_0=0$ & $f_0:=p_1+\frac{q_2^2}{4}$  \\ \hline
$\alpha_1=0$ & $f_1:=p_2$  \\ \hline
$\alpha_2=0$ & $f_2:=p_1+2p_2+t+q_1q_2-\frac{q_2^2}{4}$  \\ \hline
\end{tabular}
\end{center}
\end{proposition}
We note that when $\alpha_1=0$, we see that the system \eqref{11} admits a particular solution $p_2=0$. The system in the variables $q_1,p_1$ and $q_2$ are given by
\begin{equation}\label{pat}
  \left\{
  \begin{aligned}
   \frac{dq_1}{dt} &=2q_1^2+4p_1+2t,\\
   \frac{dp_1}{dt} &=-4q_1p_1-2\alpha_0,\\
   \frac{dq_2}{dt} &=q_2^2+4p_1-2q_1q_2.
   \end{aligned}
  \right. 
\end{equation}
This is a Riccati extension of the second Painlev\'e system in the variables $(q_1,p_1)$. Moreover, $\alpha_0=0$, we see that the system \eqref{pat} admits a particular solution $p_1=0$. The system in the variables $q_1$ and $q_2$ are given by
\begin{equation}
  \left\{
  \begin{aligned}
   \frac{dq_1}{dt} &=2q_1^2+2t,\\
   \frac{dq_2}{dt} &=q_2^2-2q_1q_2.
   \end{aligned}
  \right. 
\end{equation}
This is a Riccati extension of Airy equation in the variable $q_1$.

When $\alpha_2=0$, after we make the birational and symplectic transformations:
\begin{equation}
x_2=q_1, \ y_2=p_1+2p_2+t+q_1q_2-\frac{q_2^2}{4}, \ z_2=q_2-2q_1, \ w_2=p_2+\frac{q_1^2}{2}-\frac{q_2^2}{8}
\end{equation}
we see that the system \eqref{11} admits a particular solution $y_2=0$.

\begin{proposition}
Let us define the following translation operators{\rm : \rm}
\begin{align}
\begin{split}
&T_1:=s_1 s_2 s_1 s_0, \quad T_2:=s_1 s_0 s_1 s_2.
\end{split}
\end{align}
These translation operators act on parameters $\alpha_i$ as follows$:$
\begin{align}
\begin{split}
T_1(\alpha_0,\alpha_1,\alpha_2)=&(\alpha_0,\alpha_1,\alpha_2)+(-1,1,0),\\
T_2(\alpha_0,\alpha_1,\alpha_2)=&(\alpha_0,\alpha_1,\alpha_2)+(0,1,-1).
\end{split}
\end{align}
\end{proposition}

\begin{theorem}
Let us consider a polynomial Hamiltonian system with Hamiltonian $H \in {\Bbb C}(t)[q_1,p_1,q_2,p_2]$. We assume that

$(A1)$ $deg(H)=3$ with respect to $q_1,p_1,q_2,p_2$.

$(A2)$ This system becomes again a polynomial Hamiltonian system in each coordinate system $r_i \ (i=0,1,2)$ {\rm : \rm}
\begin{align*}
r_0:&x_0=\frac{1}{q_1}, \ y_0=-\left(\left(p_1+\frac{q_2^2}{4}\right)q_1+\alpha_0\right)q_1, \ z_0=q_2, \ w_0=p_2+\frac{q_1q_2}{2},\\
r_1:&x_1=q_1, \ y_1=p_1, \ z_1=\frac{1}{q_2}, \ w_1=-(q_2p_2+\alpha_1)q_2, \\
r_2:&x_2=\frac{1}{q_1}, \ y_2=-\left(\left(p_1+2p_2+t+q_1q_2-\frac{q_2^2}{4}\right)q_1+\alpha_2 \right)q_1,\\
&z_2=q_2-2q_1, \ w_2=p_2+\frac{q_1^2}{2}-\frac{q_2^2}{8}.
\end{align*}
Then such a system coincides with this system \eqref{11} with the polynomial Hamiltonian  \eqref{12}.
\end{theorem}
By this theorem, we can also recover the parameter's relation \eqref{13}.

We note that the condition $(A2)$ should be read that
\begin{align*}
&r_j(K) \quad (j=0,1), \quad r_2(K-q_1)
\end{align*}
are polynomials with respect to $x_i,y_i,z_i,w_i$.

\section{$D_5^{(2)}$ system} 
We study a 4-parameter family of 2-coupled $D_3^{(2)}$-systems in dimension eight with affine Weyl group symmetry of type $D_5^{(2)}$ explicitly given by
\begin{equation}\label{21}
  \left\{
  \begin{aligned}
   \frac{dq_1}{dt} &=\frac{\partial H_{D_5^{(2)}}}{\partial p_1}=2q_1^2+p_1+2t+4p_2+3p_4,\\
   \frac{dp_1}{dt} &=-\frac{\partial H_{D_5^{(2)}}}{\partial q_1}=-4q_1p_1+2q_2p_2-2\alpha_0,\\
   \frac{dq_2}{dt} &=\frac{\partial H_{D_5^{(2)}}}{\partial p_2}=q_2^2-4p_2+4p_1-2q_1q_2,\\
   \frac{dp_2}{dt} &=-\frac{\partial H_{D_5^{(2)}}}{\partial q_2}=-2q_2p_2+2q_1p_2-\alpha_1,\\
   \frac{dq_3}{dt} &=\frac{\partial H_{D_5^{(2)}}}{\partial p_3}=q_3^2-4p_3+4p_4-2q_3q_4,\\
   \frac{dp_3}{dt} &=-\frac{\partial H_{D_5^{(2)}}}{\partial q_3}=-2q_3p_3+2q_4p_3-\alpha_3,\\
   \frac{dq_4}{dt} &=\frac{\partial H_{D_5^{(2)}}}{\partial p_4}=2q_4^2+p_4+2t+3p_1+4p_3,\\
   \frac{dp_4}{dt} &=-\frac{\partial H_{D_5^{(2)}}}{\partial q_4}=-4q_4p_4+2q_3p_3-2\alpha_4
   \end{aligned}
  \right. 
\end{equation}
with the polynomial Hamiltonian \eqref{22}.

\begin{theorem}
This system \eqref{21} admits extended affine Weyl group symmetry of type $D_5^{(2)}$ as the group of its B{\"a}cklund transformations, whose generators are explicitly given as follows{\rm : \rm} with {\it the notation} $(*):=(q_1,p_1, \ldots,q_4,p_4,,t;\alpha_0,\alpha_1,\ldots,\alpha_4),$
\begin{align}
\begin{split}
s_0:(*) \rightarrow &(q_1+\frac{4\alpha_0}{4p_1+q_2^2},p_1,q_2,p_2-\frac{2\alpha_0 q_2}{4p_1+q_2^2},q_3,p_3,q_4,p_4,t;\\
&-\alpha_0,\alpha_1+2\alpha_0,\alpha_2,\alpha_3,\alpha_4),\\
s_1:(*) \rightarrow &\left(q_1,p_1,q_2+\frac{\alpha_1}{p_2},p_2,q_3,p_3,q_4,p_4,t;\alpha_0+\alpha_1,-\alpha_1,\alpha_2+\alpha_1,\alpha_3,\alpha_4 \right),
\end{split}
\end{align}
\begin{align*}
s_2:(*) \rightarrow &(q_1+\frac{2\alpha_2}{2p_1+4p_2+4p_3+2p_4+2q_1q_2-q_2q_3+2q_3q_4+4t},\\
&p_1-\frac{2\alpha_2 q_2}{2p_1+4p_2+4p_3+2p_4+2q_1q_2-q_2q_3+2q_3q_4+4t}\\
&-\frac{4\alpha_2^2}{(2p_1+4p_2+4p_3+2p_4+2q_1q_2-q_2q_3+2q_3q_4+4t)^2},\\
&q_2+\frac{4\alpha_2}{2p_1+4p_2+4p_3+2p_4+2q_1q_2-q_2q_3+2q_3q_4+4t},\\
&p_2-\frac{\alpha_2 (2q_1-q_3)}{2p_1+4p_2+4p_3+2p_4+2q_1q_2-q_2q_3+2q_3q_4+4t},\\
&q_3+\frac{4\alpha_2}{2p_1+4p_2+4p_3+2p_4+2q_1q_2-q_2q_3+2q_3q_4+4t},\\
&p_3-\frac{\alpha_2 (2q_4-q_2)}{2p_1+4p_2+4p_3+2p_4+2q_1q_2-q_2q_3+2q_3q_4+4t},\\
&q_4+\frac{2\alpha_2}{2p_1+4p_2+4p_3+2p_4+2q_1q_2-q_2q_3+2q_3q_4+4t},\\
&p_4-\frac{2\alpha_2 q_3}{2p_1+4p_2+4p_3+2p_4+2q_1q_2-q_2q_3+2q_3q_4+4t}\\
&-\frac{4\alpha_2^2}{(2p_1+4p_2+4p_3+2p_4+2q_1q_2-q_2q_3+2q_3q_4+4t)^2},\\
&t;\alpha_0,\alpha_1+\alpha_2,-\alpha_2,\alpha_3+\alpha_2,\alpha_4),\\
s_3:(*) \rightarrow &\left(q_1,p_1,q_2,p_2,q_3,p_3+\frac{\alpha_3}{p_3},q_4,p_4,t;\alpha_0,\alpha_1,\alpha_2+\alpha_3,-\alpha_3,\alpha_4+\alpha_3 \right),\\
s_4:(*) \rightarrow &(q_1,p_1,q_2,p_2,q_3,p_3-\frac{2\alpha_4 q_3}{4p_4+q_3^2},q_4,p_4+\frac{4\alpha_4}{4p_4+q_3^2},t;\\
&\alpha_0,\alpha_1,\alpha_2,\alpha_3+2\alpha_4,-\alpha_4),\\
\pi:(*) \rightarrow &(q_4,p_4,q_3,p_3,q_2,p_2,q_1,p_1,t;\alpha_4,\alpha_3,\alpha_2,\alpha_1,\alpha_0).
\end{align*}
\end{theorem}
Since these B{\"a}cklund transformations have Lie theoretic origin, similarity reduction of a Drinfeld-Sokolov hierarchy admits such a B{\"a}cklund symmetry.

\begin{proposition}
This system has the following invariant divisors\rm{:\rm}
\begin{center}\label{invariant}
\begin{tabular}{|c|c|c|} \hline
parameter's relation & $f_i$ \\ \hline
$\alpha_0=0$ & $f_0:=p_1+\frac{q_2^2}{4}$  \\ \hline
$\alpha_1=0$ & $f_1:=p_2$  \\ \hline
$\alpha_2=0$ & $f_2:=p_2+\frac{p_1+p_4}{2}+p_3+t+\frac{q_1q_2}{2}+\frac{q_3q_4}{2}-\frac{q_2q_3}{4}$  \\ \hline
$\alpha_3=0$ & $f_3:=p_3$  \\ \hline
$\alpha_4=0$ & $f_4:=p_4+\frac{q_3^2}{4}$  \\ \hline
\end{tabular}
\end{center}
\end{proposition}

\begin{proposition}
Let us define the following translation operators{\rm : \rm}
\begin{align}
\begin{split}
&T_1:=s_1 s_2 s_3 s_4 s_3 s_2 s_1 s_0, \quad T_2:=s_1 T_1 s_1, \quad T_3:=s_2 T_2 s_2, \quad T_4:=s_3 T_3 s_3.
\end{split}
\end{align}
These translation operators act on parameters $\alpha_i$ as follows$:$
\begin{align}
\begin{split}
T_1(\alpha_0,\alpha_1,\ldots,\alpha_4)=&(\alpha_0,\alpha_1,\ldots,\alpha_4)+(-2,2,0,0,0),\\
T_2(\alpha_0,\alpha_1,\ldots,\alpha_4)=&(\alpha_0,\alpha_1,\ldots,\alpha_4)+(0,-2,2,0,0),\\
T_3(\alpha_0,\alpha_1,\ldots,\alpha_4)=&(\alpha_0,\alpha_1,\ldots,\alpha_4)+(0,0,-2,2,0),\\
T_4(\alpha_0,\alpha_1,\ldots,\alpha_4)=&(\alpha_0,\alpha_1,\ldots,\alpha_4)+(0,0,0,-2,2).
\end{split}
\end{align}
\end{proposition}

\begin{theorem}
Let us consider a polynomial Hamiltonian system with Hamiltonian $H \in {\Bbb C}(t)[q_1,p_1, \ldots,q_4,p_4]$. We assume that

$(B1)$ $deg(H)=3$ with respect to $q_1,p_1, \ldots,q_4,p_4$.

$(B2)$ This system becomes again a polynomial Hamiltonian system in each coordinate system $r_i \ (i=0,1,\ldots,4)$ {\rm : \rm}
\begin{align*}
r_0:&x_0=\frac{1}{q_1}, \ y_0=-\left(\left(p_1+\frac{q_2^2}{4}\right)q_1+\alpha_0\right)q_1, \ z_0=q_2, \ w_0=p_2+\frac{q_1q_2}{2},\\
&l_0=q_3, \ m_0=p_3, \ n_0=q_4, \ u_0=p_4,\\
r_1:&x_1=q_1, \ y_1=p_1, \ z_1=\frac{1}{q_2}, \ w_1=-(q_2p_2+\alpha_1)q_2,\\
&l_1=q_3, \ m_1=p_3, \ n_1=q_4, \ u_1=p_4,\\
r_2:&x_2=q_1-\frac{q_2}{2}, \ y_2=p_1+\frac{q_2^2}{4}, \ z_2=\frac{1}{q_2},\\
&w_2=-\left(\left(p_2+\frac{p_1+p_4}{2}+p_3+t+\frac{q_1q_2}{2}+\frac{q_3q_4}{2}-\frac{q_2q_3}{4} \right)q_2+\alpha_2 \right)q_2,\\
&l_2=q_3-q_2, \ m_2=p_3-\frac{q_2^2}{4}+\frac{q_2q_4}{2}, \ n_2=q_4-\frac{q_2}{2}, \ u_2=p_4+\frac{q_2q_3}{2}-\frac{q_2^2}{4},\\
r_3:&x_3=q_1, \ y_3=p_1, \ z_3=q_2, \ w_3=p_2, \ l_3=\frac{1}{q_3}, \ m_3=-(q_3p_3+\alpha_3)q_3,\\
&n_3=q_4, \ u_3=p_4,\\
r_4:&x_4=q_1, \ y_4=p_1, \ z_4=q_2, \ w_4=p_2, \ l_4=q_3, \ m_4=p_3+\frac{q_3q_4}{2}, \ n_4=\frac{1}{q_4},\\
&u_4=-\left(\left(p_4+\frac{q_3^2}{4} \right)q_4+\alpha_4 \right)q_4.
\end{align*}
Then such a system coincides with this system \eqref{21} with the polynomial Hamiltonian  \eqref{22}.
\end{theorem}
By this theorem, we can also recover the parameter's relation \eqref{23}.

We note that the condition $(B2)$ should be read that
\begin{align*}
&r_j(K) \quad (j=0,1,3,4), \quad r_2(K-q_2)
\end{align*}
are polynomials with respect to $x_i,y_i,z_i,w_i,l_i,m_i,n_i,u_i$.

\end{document}